\documentclass[11pt]{amsart}
\usepackage{breqn,color}
\usepackage{amsmath,amsthm,amssymb,amsfonts,mathtools,bm,bbm,array,float,mathdots,dsfont,mathrsfs,latexsym}
\usepackage{caption,subcaption,multirow,longtable,lscape,wrapfig,tikz,tikz-cd,comment}
\usepackage[colorlinks=true, pdfstartview=FitV,linkcolor=blue,citecolor=blue,urlcolor=blue]{hyperref}
\usepackage{fullpage}
\usepackage{xcolor}
\usepackage[colorinlistoftodos]{todonotes}
\usepackage{verbatim}
\allowdisplaybreaks

\usepackage[square,sort,comma,numbers]{natbib}

\usetikzlibrary{patterns}

\newtheorem {theorem}    {Theorem}[section]

\theoremstyle{definition}

\newtheorem{remark}[theorem]{Remark}

\def\y{\mathbf{y}}

\numberwithin{equation}{section}
\numberwithin{figure}{section}
\numberwithin{table}{section}
\allowdisplaybreaks
\begin{document}

\title{Learning Fricke signs from Maass Form Coefficients}

\author[J. Bieri]{Joanna Bieri}
\address{University of Redlands, Redlands, CA 92373, USA}
\email{\href{mailto:joanna_bieri@redlands.edu}{joanna_bieri@redlands.edu}}

\email{}\author[G. Butbaia]{Giorgi Butbaia}
\address{Department of Physics and Astronomy, University of New Hampshire, Durham, NH 03824,
USA}
\email{\href{mailto:giorgi.butbaia@unh.edu}{giorgi.butbaia@unh.edu}}

\author[E. Costa]{Edgar Costa}
\address{
  Department of Mathematics,
  Massachusetts Institute of Technology,
  Cambridge,
  MA 02139,
  USA
}
\email{\href{mailto:edgarc@mit.edu}{edgarc@mit.edu}}
\urladdr{\url{https://edgarcosta.org}}

\author[A. Deines]{Alyson Deines}
\address{Center for Communications Research, La Jolla, USA}
\email{\href{mailto:aly.deines@gmail.com}{aly.deines@gmail.com}}

\author[K.-H. Lee]{Kyu-Hwan Lee$^{\star}$}
\address{Department of Mathematics, University of Connecticut, Storrs, CT 06269, USA  \hfill \break \indent Korea Institute for Advanced Study, Seoul 02455, Republic of Korea}
\email{\href{mailto:khlee@math.uconn.edu}{khlee@math.uconn.edu}}

\author[D. Lowry-Duda]{David Lowry-Duda$^{\dagger}$}
\address{ICERM, Providence, RI, 02903, USA}
\email{\href{mailto:david@lowryduda.com}{david@lowryduda.com}}
\urladdr{\url{https://davidlowryduda.com}}

\author[T. Oliver]{Thomas Oliver}
\address{University of Westminster, London, UK}
\email{\href{mailto:T.Oliver@westminster.ac.uk}{T.Oliver@westminster.ac.uk}}

\author[Y. Qi]{Yidi Qi$^{\ddagger}$}
\address{Department of Physics, Northeastern University, Boston, MA, USA \hfill \break \indent NSF Institute for Artificial Intelligence and Fundamental Interactions, Cambridge, MA, USA}
\email{\href{mailto:y.qi@northeastern.edu}{y.qi@northeastern.edu}}

\author[T. Veenstra]{Tamara Veenstra}
\address{Center for Communications Research, La Jolla, USA}
\email{\href{mailto:tamarabveenstra@gmail.com}{tamarabveenstra@gmail.com}}

\date{\today}

\begin{abstract}
In this paper, we conduct a data-scientific investigation of Maass forms.
We find that averaging the Fourier coefficients of Maass forms with the same Fricke sign reveals patterns analogous to the recently discovered ``murmuration'' phenomenon, and that these patterns become more pronounced when parity is incorporated as an additional feature.
Approximately $43\%$ of the forms in our dataset have an unknown Fricke sign.
For the remaining forms, we employ Linear Discriminant Analysis (LDA) to machine learn their Fricke sign, achieving $96\%$ (resp. $94\%$) accuracy for forms with even (resp. odd) parity.
We apply the trained LDA model to forms with unknown Fricke signs to make predictions. The average values based on the predicted Fricke signs are computed and compared to those for forms with known signs to verify the reasonableness of the predictions. Additionally, a subset of these predictions is evaluated against heuristic guesses provided by Hejhal's algorithm, showing a match approximately $95\%$ of the time. We also use neural networks to obtain results comparable to those from the LDA model.
\end{abstract}

\maketitle

\section{Introduction}

The purpose of this paper is to study Maass forms using various techniques from machine learning.
This builds on previous work for $L$-functions attached to number fields and arithmetic curves \cite{HLOa,HLOb,HLOc}.
Attempts to interpret this prior work led to the discovery of so-called {\em murmurations}, which are statistical correlations between the root numbers of $L$-functions and their Dirichlet coefficients, first observed in the context of elliptic curves~\cite{HLOP}.
Murmurations for Maass forms were established in~\cite{BLLDSHZ}.

We show that supervised machine learning techniques can be used to predict the Fricke sign of a Maass form based on its coefficients, \emph{without} indicating the level. In contrast to the papers cited previously, we note that unsupervised techniques such as $k$-means clustering and PCA were unsuccessful at producing clusters separated by Fricke sign.
The relationship between murmurations of Maass forms and machine learning their Fricke sign is analogous to the relationship between murmurations of elliptic curves and machine learning their ranks (as in~\cite{HLOc,HLOP}).

The machine learning experiments presented here use the database of
Maass forms from the LMFDB~\cite{lmfdb}.
Among the 35,416 Maass forms in the LMFDB, some 15,423 of them lack rigorously computed Fricke signs, though it is possible to guess them heuristically.
The computations of Maass forms underpinning the LMFDB rely on a combination of different techniques~\cite{dld2025database}, including automorphic certification~\cite{C22}, explicit trace formulas~\cite{seymourhowell}, and an implementation of Hejhal's algorithm~\cite{hejhal} described in forthcoming work~\cite{LDSH}.
After that, the Fricke sign is subsequently determined by the signs of the coefficients $a_p$ for each prime $p$ dividing the level of the Maass form.
In practice, however, it is often quite hard to compute the Maass form eigenvalue and coefficients with sufficient precision to directly recover the Fricke sign. This difficulty explains why a substantial portion of the dataset currently has unknown Fricke signs.

Consequently, it is
surprising that our machine learning results achieve high accuracy ($>96\%$) in predicting the Fricke sign.
Our approach fundamentally differs from conventional methods:
we consider the dataset collectively and employ machine learning tools instead of relying solely on the intrinsic properties of individual forms. Moreover, we apply the trained machine learning model to Maass forms with unknown Fricke signs in order to predict them.
By comparing these predictions with heuristic guesses, we demonstrate that our results are
reasonable and may provide valuable information that could help determine the signs precisely.

Furthermore, the decomposition of the Fricke sign into a product of local factors implies that it may be viewed as a \textsl{parity function} in the sense of \cite{SS25}, where it is noted that general purpose machine learning methods have not previously succeeded in learning such functions from sign patterns.
In Section \ref{sect:no_cheat}, we will show that, although sign patterns are implicit in our data presentation, our classifiers are learning something more from the input features.

In future work, it will be beneficial to interpret the machine learning results, to identify which features contribute most significantly to the high accuracy, and to understand how they do so.
This analysis may lead to a deeper understanding of the Fricke sign. Additionally, we anticipate that other important quantities in number theory can also be effectively approached through machine learning. We hope this methodology will open new avenues for studying various objects that are difficult to compute in conventional ways.

We conclude this introduction with a summary of each section.
In Section~\ref{s.Maassforms}, we review the necessary definitions.
In Section~\ref{s:averaging}, we observe murmurations of Maass forms by averaging their coefficients.
In Section~\ref{s:LDA}, we use Linear Discriminant Analysis (LDA) to predict the Fricke sign based upon these same coefficients.
In Section~\ref{sect:no_cheat}, we outline a strategy for predicting the Fricke sign of a Maass form based upon factorising its level, and confirm that this approach is not what is taking place in the LDA.
In Section~\ref{s:NN}, we train a neural network on the coefficients and the spectral parameter and achieve an accuracy comparable to that attained with LDA.
In Section~\ref{s:comparisonLMFDB}, we compare our predictions to those based on a heuristic version of Hejhal's algorithm.

\section*{Acknowledgments}

Costa was supported by Simons Foundation grants 550033 and SFI-MPS-Infrastructure-00008651.
Lee was supported by Simons Foundation grant 712100.
Lowry-Duda was supported by Simons Foundation grant 546235.
Qi was supported by the NSF grant PHY-2019786 (the NSF AI Institute for Artificial Intelligence and Fundamental Interactions).

We would also like to express our gratitude to the organizers of the Mathematics and Machine Learning Program at CMSA during the Fall 2024 semester, where this project was initiated.

\section{Experiments with Maass forms}

\subsection{Definitions}\label{s.Maassforms}
For a broad overview, the reader is referred to~\cite{FL05}; a more extensive (but less approachable) reference is~\cite{DFI}.
Let $\Delta$ denote the Laplace--Beltrami operator on the upper half-plane.
A weight $0$ Maass cuspform $f$ on $\Gamma_0(N)$ is a smooth square-integrable function $f:\mathcal{H}\rightarrow\mathbb{C}$ satisfying $f(\gamma z)=f(z)$, for all $\gamma\in\Gamma_0(N)$, and $(\Delta-\lambda)f(z)=0$, for some $\lambda\in\mathbb{C}$.
We refer to $N$ as the \textsl{level} of $f$.
Assuming the Selberg eigenvalue conjecture, we may write $\lambda=\frac14+R^2$ for $R\in\mathbb{R}_{\ge 0}$.
In what follows, we refer to $R$ as the \textsl{spectral parameter}.

If $f$ is a Maass form on $\Gamma_0(M)$ and $M$ divides $N$, then, for any $k$ dividing $N/M$, $f(kz)$ is a so-called \textsl{oldform} on $\Gamma_0(N)$.
We will consider only Maass forms that are not oldforms, known as \textsl{newforms}, which live naturally on $\Gamma_0(N)$.
We write $S_0^{\textup{new}}(\Gamma_0(N))$ for the set of weight $0$ Maass newforms on $\Gamma_0(N)$.

It is known that $f\in S_0^{\textup{new}}(\Gamma_0(N))$ has a Fourier expansion of the form
\begin{equation}\label{eq.maass_exp}
  f(x+iy) = \sum_{n\neq0}a_n\sqrt{y}K_{iR}(2\pi|ny|)\exp(2\pi inx),
\end{equation}
where $K_{iR}(u)=\frac12\int_0^{\infty}\exp(-|u|(t+t^{-1})/2)t^{iR-1}dt$ is a modified Bessel function of the second kind.
The Maass form $f$ satisfies the functional equation
\[f(z)=w_Nf\left(-\frac{1}{Nz}\right),\]
for $w_N\in\{\pm1\}$.
We refer to $w_N$ as the \textsl{Fricke sign}.

There is an involution on $\mathcal{H}$ given by reflection in the imaginary axis, $z\mapsto-\bar{z}$.
We call a Maass form even (resp.\ odd) if  $f(-\bar{z})=f(z)$ (resp. $-f(z)$).
Let $\sigma(f)$ be $0$ (resp. $1$) if $f$ even (resp. odd).
Applying equation~\eqref{eq.maass_exp} for a Maass form $f$ with fixed parity, we deduce:
\begin{equation}\label{eq.Maass-expansion}
  f(x+iy) = \sum_{n=1}^{\infty}a_n\sqrt{y}K_{iR}(2\pi \lvert ny \rvert)
  \begin{cases}
    2\cos(2\pi inx),    &\sigma(f)=0,\\
    2i\sin(2\pi inx),   &\sigma(f)=1.
  \end{cases}
\end{equation}
The $L$-function associated to $f$ is given by $L_f(s)=\sum_{n=1}^{\infty}a_nn^{-s}$, which converges for $\mathrm{Re}(s) > 1$.
The completed $L$-function is:
\[
    \Lambda_f(s)
    =
    \left(\frac{\sqrt{N}}{\pi}\right)^s
    \Gamma\left(\frac{s+\sigma(f)+iR}{2}\right)
    \Gamma\left(\frac{s+\sigma(f)-iR}{2}\right)
    L_f(s).
\]
The completed $L$-function satisfies the functional equation
\[\Lambda_f(s)= \epsilon \bar{\Lambda}_f(1-s),\]
where $\epsilon$ is the root number.
In particular, we have $\epsilon = (-1)^{\sigma(f)} w_N$ (see \cite[\S8]{DFI}).

Given complete information about the coefficients, the Fricke sign is easily computable;
on $\Gamma_0(N)$ with $N$ squarefree, the coefficient $a_N$ encodes the Fricke sign. However, explicitly computing Maass form coefficients is extremely hard in practice.
Conjecturally, all data associated to a generic Maass form is transcendental and independent of typical transcendental constants (cf.~\cite{BSV}).

The Fricke signs
in the LMFDB were rigorously computed from the coefficients, and are missing when these coefficients were not computed with sufficient precision.
The challenge most often starts with numerical difficulties in
the explicit trace formula~\cite{seymourhowell}.
Increasing the precision of the trace formula requires computing many
class numbers rigorously, which is
very hard.
A different approach would be to implement the confirmation algorithm from~\cite{C22} for general level.
Whilst there are other methods to determine the Fricke sign of a Maass form, none have been implemented rigorously.

In the following subsections, we will show that one can use methodology from machine learning to predict the Fricke sign from finitely many coefficients of a Maass newform.
This is analogous to machine learning the rank (parity) of an elliptic curve, which is also manifest in the sign of a completed $L$-function \cite{HLOc,HLOP}.
Furthermore, this is connected to murmurations of Maass forms as in \cite{BLLDSHZ}, which amount to subtle statistical correlations between the sign and the coefficients.
Finally, we will compare our predictions with the heuristic results from Hejhal's algorithm described just above.

\subsection{Averaging}\label{s:averaging}

Our analysis uses a dataset $\mathcal{L}$ containing the 35,416 rigorously computed Maass forms from the LMFDB \cite{lmfdb}.
The dataset contains the first 1,000 Fourier coefficients $a_n$ for every Maass form, each of which has weight $0$, trivial character, and integral level $N$ in the range from $1$ to $105$ \cite{maassdata}.
In the analysis that follows, we are particularly interested in the parity $\sigma$  and the Fricke sign $w$ (we will often drop the subscript from $w_N$ if there is no need to specify $N$).
As discussed above, in the LMFDB, the value of the Fricke sign is not always rigorously computed, in which case it is given the value $w=0$.
The breakdown of the Maass forms dataset $\mathcal L$ with different values for $\sigma$ and $w$ are given in Table \ref{tab:maass_counts}.

\begin{table}
\begin{center}
\begin{tabular}{|c|c|c|c|c|}
\hline
& $w=-1$ & $w=1$ & $w=0$ (unknown) & Total \\
\hline
$\sigma(f)=0$ (even) &$5009$&$7171$&$6173$&$18353$\\
\hline
$\sigma(f)=1$ (odd)  &$3724$&$4089$&$9250$&$17063$ \\
\hline
Total & $8733$ & $11260$ & $15423$ & $35416$\\
\hline
\end{tabular}
\caption{Counts of Maass forms by parity and Fricke sign.}
\label{tab:maass_counts}
\end{center}
\end{table}

\begin{figure}[h]
\centering
\includegraphics[width=.45\textwidth]{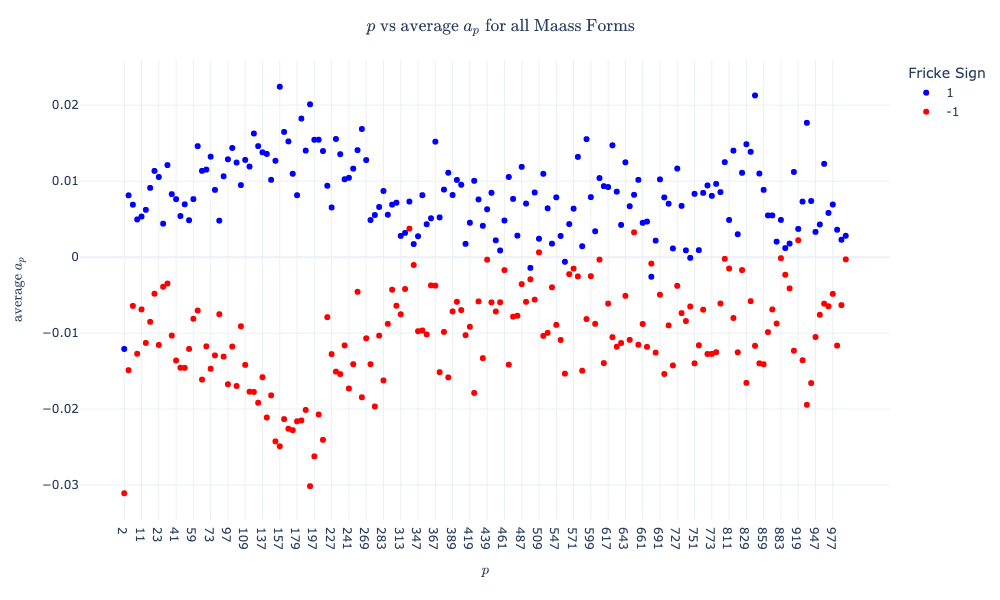}
\includegraphics[width=.45\textwidth]{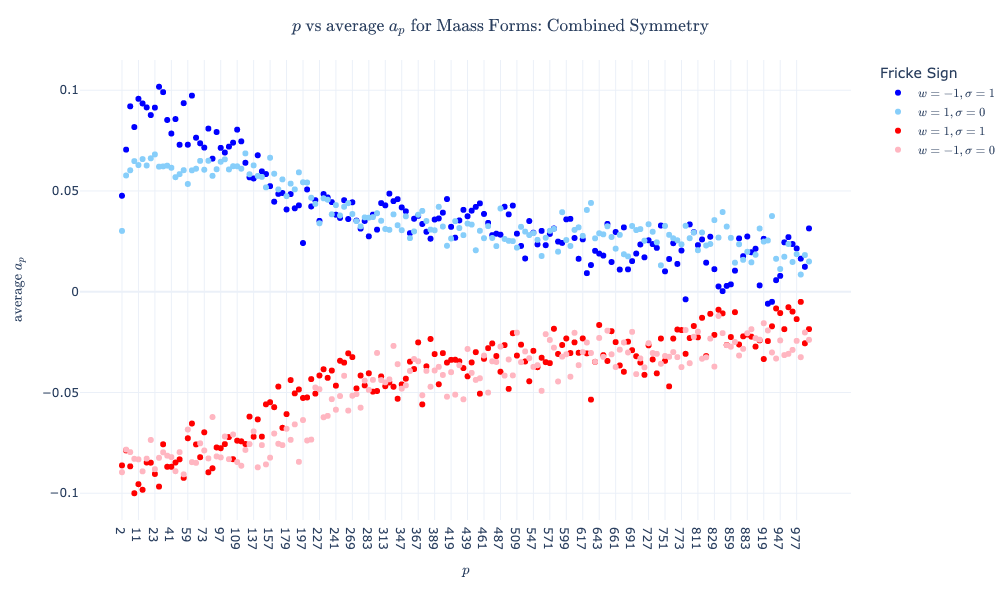}
\caption{Average value of $a_p$ over Maass forms with given Fricke sign, with and without separating by symmetry}
\label{fig:allsym-nonorm}
\end{figure}

\begin{figure}[h]
\centering
\includegraphics[width=.8\textwidth]{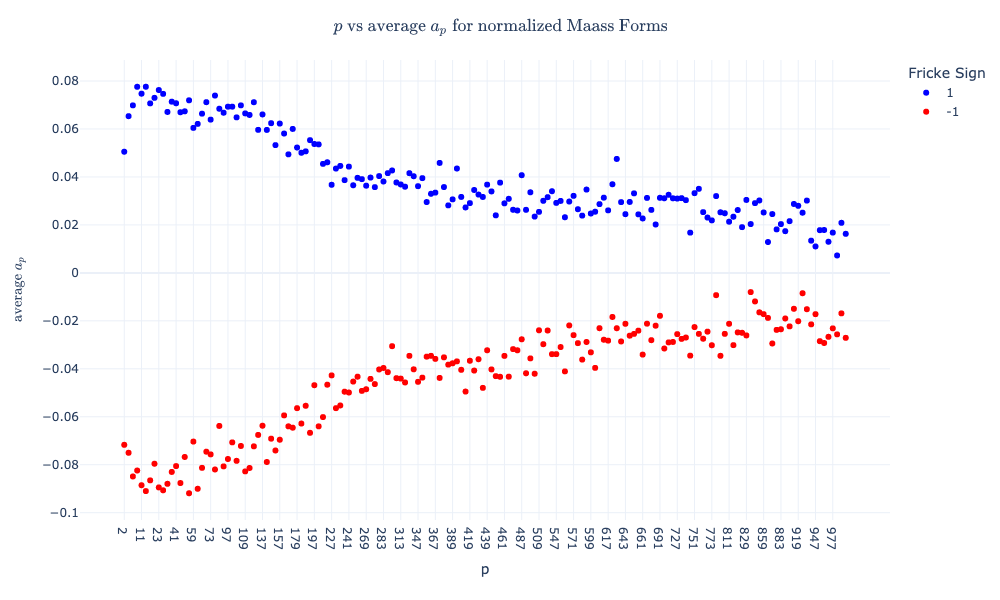}
\caption{Average value of $(-1)^{\sigma(f)} a_p$ over Maass Forms separated by Fricke sign.}
\label{fig:allsym}
\end{figure}

We may write $\mathcal{L}$ as a disjoint union $\mathcal{L}_0\coprod\mathcal{L}_1\coprod\mathcal{L}_{-1}$, in which $\mathcal{L}_i=\{f\in\mathcal{L}:w=i\}$.
Figures \ref{fig:allsym-nonorm} and \ref{fig:allsym} provide clear evidence of murmuration-like distinction between the forms in $\mathcal{L}_1$ and those in $\mathcal{L}_{-1}$.
More precisely, in Figure \ref{fig:allsym-nonorm},  for primes $p<1000$, we plot the average value of $a_p$  over $\mathcal{L}_1$ and $\mathcal{L}_{-1}$, where $a_p$ is as in equation~\eqref{eq.Maass-expansion}.  We note that the separation is much better when also taking symmetry into account. This is equivalent to the idea of separating by root
number. Because the root number $\epsilon = (-1)^{\sigma(f)} w$ we explore normalizing the coefficients
by multiplying by
$(-1)^{\sigma(f)}$.
In Figure~\ref{fig:allsym} we see this produces an even more
 distinctive murmuration pattern.

\subsection{Predicting the Fricke sign with LDA}\label{s:LDA}

The clear separation in Fricke sign provided by averaging $(-1)^{\sigma(f)}a_p$ indicates that it may be possible to predict the Fricke sign based on these features, perhaps even using a fairly simple technique.
We choose to undertake Linear Discriminant Analysis (LDA) to learn a linear decision boundary between classes of Fricke sign.
LDA is a supervised dimensionality reduction and classification technique based on Bayes theorem \cite[Section~4.3]{Hastie}.

By writing each $f\in\mathcal{L}$ as in equation~\eqref{eq.Maass-expansion}, we may construct the following $1000$-dimensional feature vector:
\[\mathcal{D}=\left\{\left((-1)^{\sigma(f)}a_n\right)_{i=1}^{1000}:f\in\mathcal{L}\right\}\subset\mathbb{R}^{1000}.\]
We may write $\mathcal{D}=\mathcal{D}_0\coprod\mathcal{D}_1\coprod\mathcal{D}_{-1}$, where $\mathcal{D}_i$ contains only vectors corresponding to forms in $\mathcal{L}_i$.

\begin{figure}
    \centering
    \includegraphics[width=.45\textwidth]{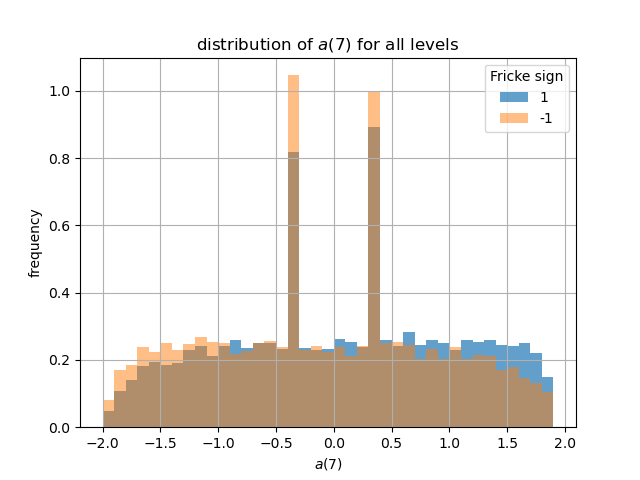}
\includegraphics[width=.45\textwidth]{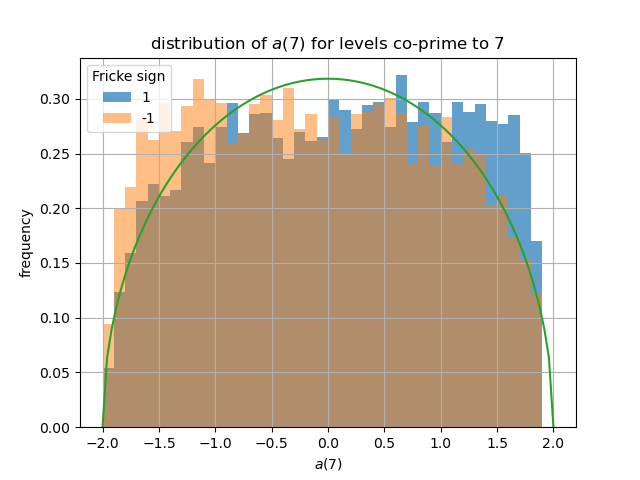}
    \caption{Comparing the distribution of $a_7$ for Maass forms in $\mathcal{L}_1$ and $\mathcal{L}_{-1}$.
    In the left (resp. right) frame, we consider all levels (resp. only levels co-prime to $7$).
    The brown represents areas where the histograms overlap, and the green arc in the right frame is the semi-circle $
y=\frac{1}{2\pi}\sqrt{4-x^2}$, given by the (vertical) Sato--Tate distribution \cite{Sar}.}
    \label{fig:fricke_dist_7}
\end{figure}

LDA is likely to work well when the covariance is the same for all classes.
To determine if LDA is a good technique for $\mathcal{D}$, we examined the covariance of $\mathcal{D}_1$ and $\mathcal{D}_{-1}$.
Visual inspection of covariance matrices, their eigenvalues, and their determinants exhibited similar results between the two classes.
To be more rigorous, we also applied Box's M test~\cite{box49mtest} for each feature $a_n$.
Without any further feature engineering, Box's M test declared equal covariance for $641$ of the $1000$ features, including all prime indices.
We note that the distribution of the features $a_p$ indexed by primes $p$ has spikes because  the fixed value $a_p = -w_p/\sqrt{p}$ is repeated when $p$ divides the level (cf. Section~\ref{sect:no_cheat}).
By way of example, in Figure \ref{fig:fricke_dist_7}, we plot the distribution of $a_p$ in the case that $p=7$, firstly over all forms in $\mathcal{L}_{\pm1}$ and then only over those with level coprime to $7$.
With this in mind, for each $n$, we removed from $\mathcal{L}$ all Maass forms for which $n$ divides its level.
Upon doing so, Box's M test declared that the two classes had equal covariance for all but $33$ of the $1000$ features, all of which had composite index.
Altogether, we are satisfied that the classes
have equal covariance for the vast majority of the $n$.
We further explore the impact of the values of $a_n$ for $n$ dividing the level in Section \ref{sect:no_cheat}.

We ran supervised learning experiments on the labeled dataset $\mathcal{D}_1\coprod\mathcal{D}_{-1}$.
In each experiment, we split a subset into testing, training, and validation sets, and the testing set accounted for $20\%$ of the relevant data with the training-validation split also 80-20.
Firstly, we trained the LDA using all available training data ($12795$ observations) and recorded $96.1\%$ accuracy on the validation data.
Secondly, we masked the training data so as to include only even forms ($7772$ observations) and recorded $94.9\%$ accuracy on the similarly masked validation data.
Thirdly, we trained only data from odd forms ($5023$ observations), and recorded $96.3\%$ accuracy on the similarly masked validation data.
In all three cases, the accuracy on the testing set is almost identical to that on the validation set.
These experiments indicate that, without any hyperparameter tuning, we are able to get high predictive accuracy.

\begin{figure}[h]
\centering
\includegraphics[width=.8\textwidth]{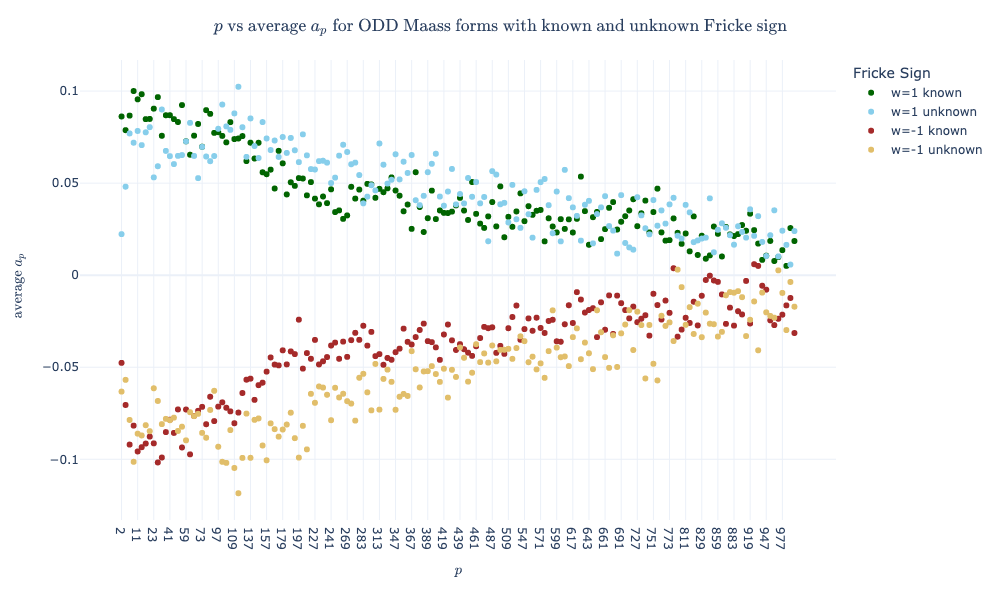}
\caption{Average value of $a_p(-1)^{\sigma(f)}$ for Maass forms with odd parity, separated by Fricke sign. }
\label{fig:LDApredictions}
\end{figure}

\begin{figure}[h]
\centering
\includegraphics[width=.8\textwidth]{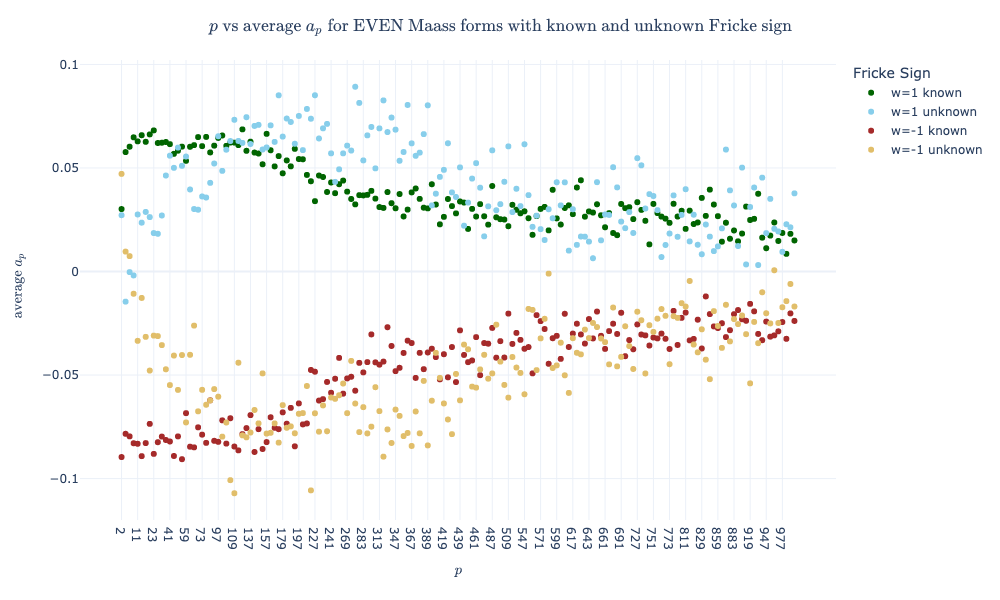}
\caption{Average value of $a_p(-1)^{\sigma(f)}$ for Maass forms with even parity, separated by Fricke sign. }
\label{fig:LDApredictions_even}
\end{figure}

Building on the paragraph above, we also make predictions using the trained LDA model for the Maass forms with unknown Fricke sign, that is, those in $\mathcal L_0$.
Then, in order to check whether our predictions are reasonable,
we examined if the average value of $(-1)^{\sigma(f)}a_p$ for the newly predicted Fricke sign follows similar murmuration patterns as the rigorously proved (known) values.
Figures \ref{fig:LDApredictions} and \ref{fig:LDApredictions_even} demonstrate that our predicted values mimic the murmuration patterns seen in Figure~\ref{fig:allsym-nonorm}.
By comparing Figures \ref{fig:LDApredictions} and \ref{fig:LDApredictions_even} we make the interesting observation that, when $p$ is small, there is greater similarity between genuine and predicted murmurations for odd Maass forms than there is for even Maass forms.

\subsection{Embedding the Fricke signs into the Fourier coefficients}\label{sect:no_cheat}
The (global) Fricke sign introduced in Section~\ref{s.Maassforms} may be written as a product of the form
\begin{equation}\label{eqn:global_fricke}
 w_N=\prod_{p|N} w_p,
\end{equation}
in which, for a prime number $p$, the number $w_p$ is the local Fricke sign.
If $p$ divides the level $N$ then we have
\begin{equation}\label{eqn:local_fricke}
a_p=\frac{-w_p}{\sqrt{p}}.
\end{equation}
We conclude from equation~\eqref{eqn:local_fricke} that the Fricke sign is determined by the sequence $(\mathrm{sgn}(a_p))$ of signs and the number of prime factors of $N$.

If the Fricke sign is unknown, we do not know the value of
$a_p$ for $p|N$, and these values are defined to be $0$ in the LMFDB.
Similarly, when $\gcd(n,N)>1$, we have $a_n=0$.
Given that we are training on a dataset where the Fricke sign is known, one might wonder whether perhaps LDA is just learning how to determine the Fricke sign from the $a_n$ with $\gcd(n,N)>1$.
To clarify, we do experiments that show this is not the case.
Before describing the experiments, we note that, for level $1$ Maass forms (of which there are $2202$ in our dataset), as $\gcd(n,1)=1$, the Fourier coefficients do not directly encode information about the Fricke sign in the above way.
Similarly, as the largest level in our dataset is $105$, for any of the primes $p>105$, the value of $a_p$ does not directly encode this information.

In our previous section we used vectors of the form $(a_n)_{n=1}^{1000}$ to make predictions with LDA.
In this section, we will explore, firstly, the impact of only using coefficients with prime index, and, secondly, the impact of removing the extra information about the Fricke sign embedded in the coefficients when $\gcd(n,N)>1$ by setting all of these coefficients to zero.
Subsequently, we examine the accuracy of both the first and second experiments based on the number of coefficients used to train the data.

\begin{figure}[h]
\centering
\includegraphics[width=.8\textwidth]{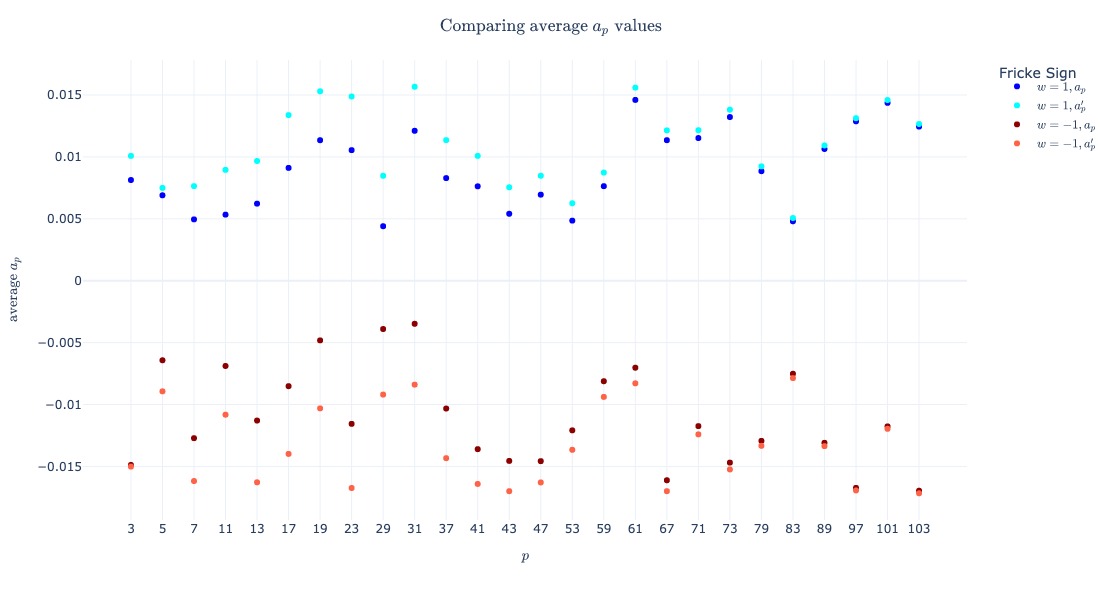}
\caption{Comparing average values of $a_p
$ versus $a'_p$ by Fricke sign}
\label{fig:levelanalysis}
\end{figure}

In Figure \ref{fig:levelanalysis}, we compare the average value of $(-1)^{\sigma(f)}a_p$ and $(-1)^{\sigma(f)}a'_p$ over $f\in\mathcal{L}_1$ and $f\in\mathcal{L}_{-1}$, where
\[a'_n=\begin{cases}0,&\gcd(n,N)>1,\\
a_n,&\gcd(n,N)=1.\end{cases}\]
This matches the data that we are trying to predict, as all $a_n=0$ if $\gcd(n,N)\neq 1$ for the Maass forms with unknown Fricke sign.
For $p>105$ we have $a_p=a'_p$, and so we restrict the graph to the relevant primes.
We see that there is very little impact to the average, and still good separation between the different Fricke signs.

In Table \ref{tab:no_cheat}, we record the accuracy of LDA when applied to vectors of the following form:
\[(a_n)_{n=1}^{1000},~(a'_n)_{n=1}^{1000},~(a_{p})_{p < 1000},~(a'_{p})_{p < 1000},\]
where the indices range over the $168$ primes less than $1000$.
In these experiments, we see that replacing certain coefficients by zero does not make much difference to the accuracy.
This is good news; we believe this means LDA should do well at predicting the Fricke sign in the case where it is unknown.
Since the values of the coefficients are multiplicative, we would expect the $a_p$ to suffice.
It is therefore surprising that there is a substantial improvement in our accuracy when using all $a_n$ rather than just $a_p$.
We will revisit this theme below in greater depth.

\begin{table}
\begin{tabular}{|c|c|c|c|}
\hline
Features & all parity (normalized) & even parity  & odd parity  \\
\hline
$a_n$ & 0.9612  & 0.9488 & 0.9633\\
\hline
$a'_n$  & 0.9456  & 0.9564 & 0.9633\\
\hline
$a_{p_i}$ &
0.8625 &
0.8261 &
0.8800\\
\hline
$a'_{p_i}$ &
0.8615 &
0.8329 &
0.8769 \\
\hline
\end{tabular}
\caption{Comparing the accuracy of LDA predictions for different features and parities.}
\label{tab:no_cheat}
\end{table}

We next analyze how these methods compare with each other when predicting the Fricke sign of Maass forms in $\mathcal{L}_0$.
Predictions based on $(a_n)_{n=1}^{1000}$ and $(a_{p})_{p < 1000}$ agree $83.29\%$ of the time.
Predictions based on $(a_n)_{n=1}^{1000}$ and $(a'_n)_{n=1}^{1000}$ agree $97.50\%$ of the time.
Predictions based on $(a_{p})_{p < 1000}$ and $(a'_{p})_{p < 1000}$ agree $99.81\%$ of the time.
This is consistent with the accuracy results outlined above in that  zeroing out the coefficients $a_n$ when $\gcd(n,N)>1$ does not have much impact on the predictions, but using $a_n$ versus $a_p$ does have a substantial impact on the predictions.

\begin{figure}[h]
\centering
\includegraphics[width=.8\textwidth]{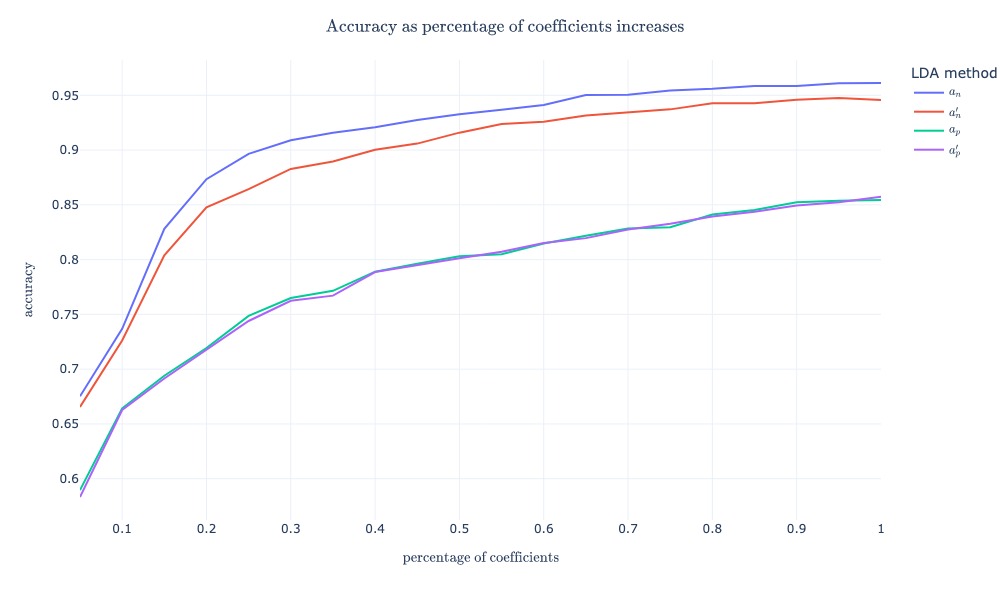}
\caption{Accuracy on validation set as number of $a_n$ increases}
\label{fig:numap_nozeros}
\end{figure}

Finally, we analyze the different LDA methods used in this section and determine their accuracy relative to the number of coefficients used to train the data; see Figure \ref{fig:numap_nozeros}.
We again notice that LDA using $(a_n)_{n=1}^{1000}$ has better accuracy than just using $(a_{p})_{p < 1000}$, and that there is little difference across the board for correct coefficients versus zeroed out coefficients.
Interestingly, we note that there is striking improvement in accuracy for including more $a_n$ with small $n$, but then a leveling off as $n$ increases.
This still occurs but is less pronounced when using $a_p$.

As mentioned earlier, it is surprising that LDA performs better when using $(a_n)_{n=1}^{1000}$ than when using  just $(a_{p})_{p < 1000}$.
Indeed, the Fourier coefficients are multiplicative, and so there isn't more information in the composite indices.
We ran several tests to see if we could explain this observation.
One might speculate that $a_p$ for small primes $p$ are likely to have a larger impact on classification.
If that were the case, then perhaps the LDA performance improved with the inclusion of composite indices because the multiplicative property was accentuating the power of prediction coming from  $a_p$ with small prime factors.
To investigate this, we represented each Maass form by vectors of the form $(a_n)_{n\in S}$, where
\[S=\{ n \leq 1000 \; | \; n \text{ is divisible by } 2,3 \text{ or } 5 \}.\]
Unexpectedly, LDA performed poorly with this presentation, with only $78.5\%$ accuracy, even though it included $733$ of the $1000$ features in $\mathcal{D}$.
Also including all the primes in this range, that is, using vectors of the form $(a_n)_{n\in S'}$ where
\[S'=\{ n \leq 1000 \; | \; n \text{ is divisible by } 2,3 \text{ or } 5 \text{ or $n$ is prime} \},\]
the accuracy increased to
$89.8\%$ with $898$ of the $1000$ original features.
While this is a mild improvement over the $86.2\%$ accuracy when using only prime indices, it is nowhere near the $96.1\%$ accuracy when using all indices $\leq1000$, even though it contains $90\%$ of those values.
It seems clear that extra information from $a_p$ with small primes $p$ is not improving the performance of LDA.
Perhaps, LDA is learning something useful from the multiplicative property itself.
The subset of indices we found with the best accuracy was
\[S''=\{ n \leq 1000 \; | \; n \text{ is divisible by 1 or 2 prime factors}\}.\]
For this presentation, LDA predicted Fricke signs with an accuracy of $95.3\%$ using $702$ of the 1000 features.
A few other interesting results are shown in Table \ref{tab:coeff_subset}.

\begin{table}[h]
    \centering
    \begin{tabular}{|c|c|}
    \hline
     Feature indices & LDA accuracy \\
     \hline
     $\{n\in\mathbb{Z}:1\leq n\leq1000\}$ &  96.1\% \\
    $\{n\text{ prime}:1\leq n\leq1000\}$ &  86.2\% \\
    $\{n\text{ $45$-smooth}:1\leq n\leq1000\}$&  75.3\% \\
    $\{n\text{ even}:1\leq n\leq1000\}$ &  70.6\% \\
    $\{n\text{ odd}:1\leq n\leq1000\}$ &  93.4\% \\
    $\{n\in\mathbb{Z}:1\leq n\leq500\}$ &  93.3\% \\
    \hline
    \end{tabular}
    \caption{Accuracy of LDA when predicting Fricke sign from features indexed by various subsets of integers.
    All of the subsets in the bottom portion of the table have about $500$ elements. In particular, smoothness bound of $45$ was chosen with that in mind.}
    \label{tab:coeff_subset}
\end{table}

\subsection{Other LDA analysis}
\begin{figure}[htb]
    \centering
    \includegraphics[width=.45\textwidth]{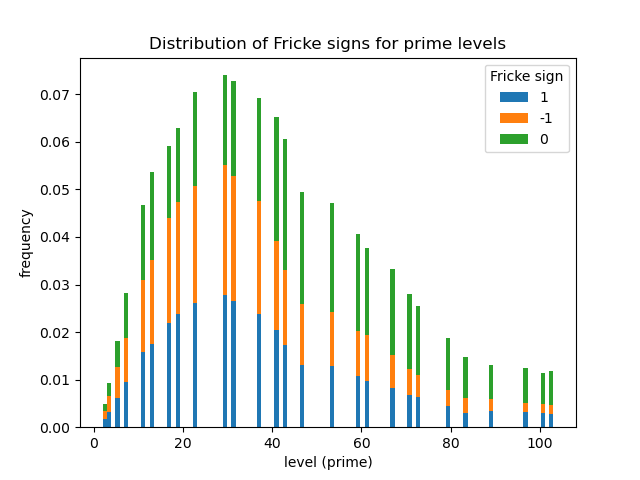}
\includegraphics[width=.45\textwidth]{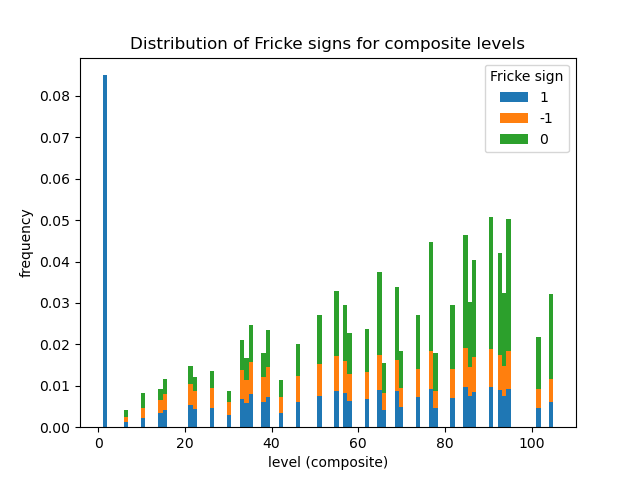}
    \caption{Comparing the distribution of Fricke signs by level. Unknown Fricke signs are denoted by 0.}
    \label{fig:fricke_signs_by_level}
\end{figure}

It is possible that the success of LDA could be explained by some subset of Maass forms that are especially easy to predict. For example, perhaps Maass forms of some levels are easier to predict than others.
It is easy to check that the distribution of Fricke signs by level is generally even between $\mathcal{L}_1$ and $\mathcal{L}_{-1}$ (cf. Figure~\ref{fig:fricke_signs_by_level}).
The only exception is level $N=1$, where the Fricke sign is always $1$.

\begin{figure}[htb]
    \centering
    \includegraphics[width=0.5\linewidth]{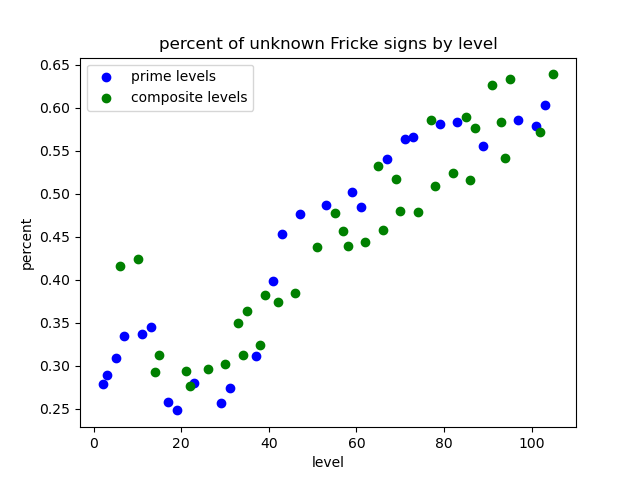}
    \caption{Comparing the percentages of known and unknown Fricke sign by level.
    }
    \label{fig:unk_v_known}
\end{figure}

We might expect forms on prime level to be easier to predict than forms on composite level, as the global Fricke sign is a product of local Fricke signs corresponding to each prime factor of the level.
Further, the overall quality of the Maass forms approximations in the LMFDB is better for lower level than higher level.
(This is caused by multiple compounding factors, including less separation between consecutive eigenvalues, forcing Heisenberg-Uncertainty tradeoffs and less well-behaved test functions in the trace formula used to compute the Maass forms).
This leads the the general positive slope in Figure~\ref{fig:unk_v_known}.


\begin{figure}[htb]
    \centering
    \includegraphics[width=0.7\linewidth]{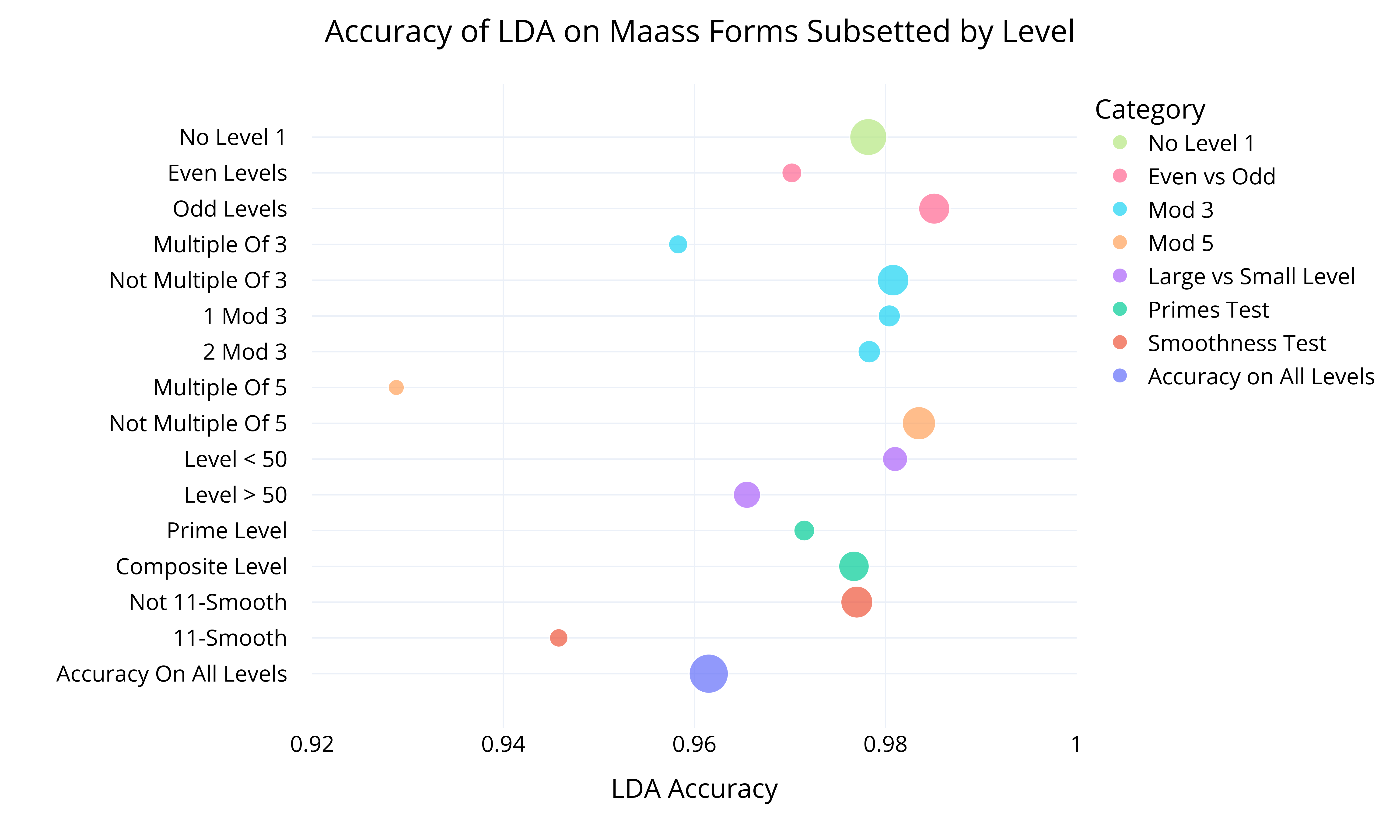}
    \caption{Exploring accuracy of LDA on subsets of Maass forms for given levels. The size of the dot corresponds to the size of the training set. }
    \label{fig:level_subset}
\end{figure}

Building on this theme, we experimented with many different subsets of Maass forms of various levels and did not find any that were especially easy or hard to predict.
The only exception was that if the subset was really small then it was, not surprisingly, hard to predict as there was not enough training data.
These experiments used all the features and are summarized in Figure \ref{fig:level_subset}.

\subsection{Predicting the Fricke sign with a neural network}\label{s:NN}
\begin{figure}[htb]
    \centering
    \includegraphics[width=0.8\textwidth]{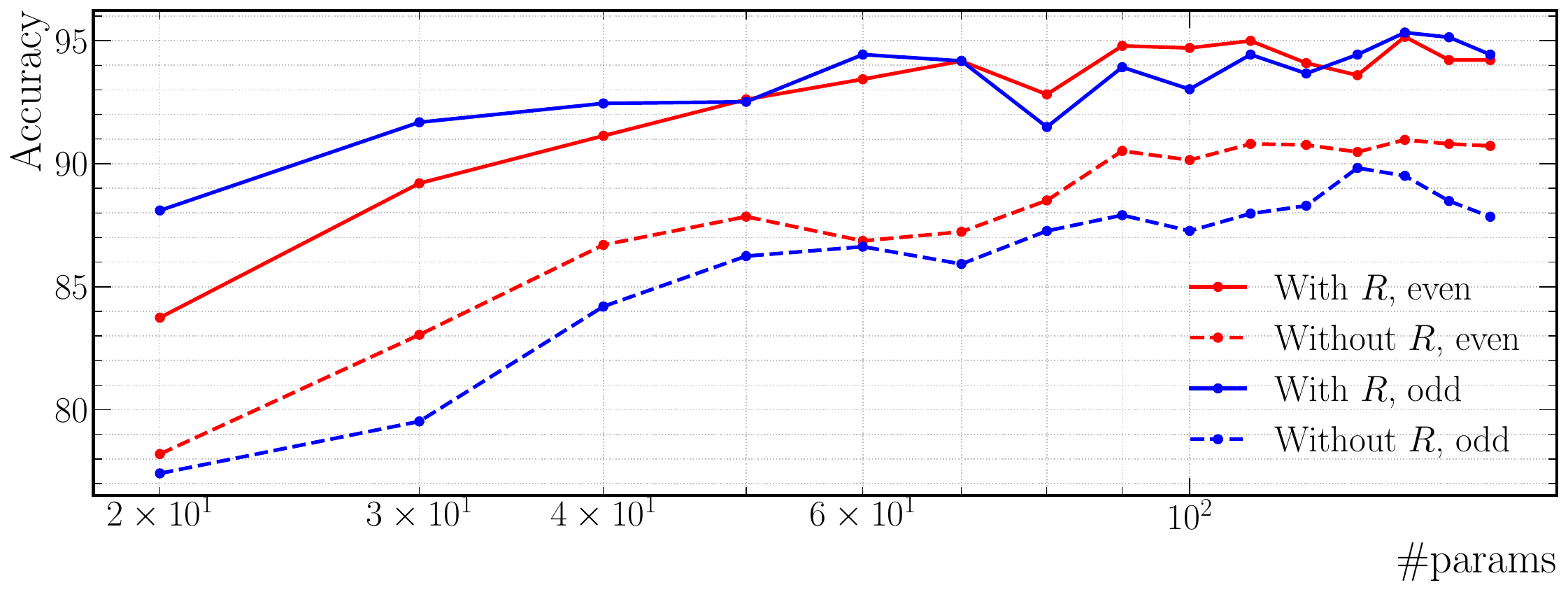}
    \caption{Dependence of the accuracy on the number of the coefficients $a_p$ (for $p$ prime) used as the input to the neural network trained for predicting Fricke sign. Here $R$ denotes the spectral parameter.}
    \label{fig:acc_vs_params_NN}
\end{figure}

Motivated by the success of LDA, we trained neural networks to predict the Fricke sign
using feature vectors of the form $(a_2,a_3,\dots,a_{p_{d}},R)$.
As in the previous experiments, we take $d=168$, so as to use all primes $<1000$.
In fact, we construct two different neural networks, one for each possible parity of the form.
This is motivated by
Figure.~\ref{fig:allsym-nonorm}.
Unlike in the previous sections, we have included the spectral parameter $R$ as a feature.
Since the analytic conductor of a Maass form is (roughly) equal to $\frac{NR^2}{4\pi^2}$, experiments using $R$ as a feature alongside the coefficients are similar to the experiments in \cite{KV22} (in which the elliptic curve conductor is used alongside the Frobenius traces to predict the rank).
We note that neglecting $R$ significantly reduced the accuracy of the predictions (cf. Figure \ref{fig:acc_vs_params_NN}).

\begin{figure}
	\centering
	\begin{tikzpicture}[shorten >= 1pt, ->, draw=black, line width=1.5]
		\tikzstyle{RR}=[very thick, rounded corners, draw=black, minimum size=25];
		\tikzstyle{Split} = [rectangle split, rectangle split parts=2, rectangle split part fill={blue!30,red!20}, rounded corners, draw=black, very thick, minimum size=23, inner sep=5pt, text centered];
		\tikzstyle{neuron}=[very thick,circle,draw=black,minimum size=25,inner sep=0.5,outer sep=0.6]
		\foreach \name / \y in {1,...,4}
			\node[RR] (I-\name) at (0,-1*\y cm) {
				\ifnum \y = 1
					$a_2$
				\fi
				\ifnum \y = 2
					$a_3$
				\fi
				\ifnum \y = 3
					$...$
				\fi
				\ifnum \y = 4
					$a_{p_d}$
				\fi
				\ifnum \y = 5
					$R$
				\fi
			};

		\node[RR, pattern=north west lines, pattern color=red!30] (I-5) at (0,-1*5 cm) {$R$};

		\foreach \name / \y in {1,...,5}
			\path[yshift=(1)*0.5cm]
				node[neuron] (N1-\name) at (3, -1.2*\y cm) {
					{\small $\mathsf{n}_\y^{(1)}$}
         	};

		\foreach \i in {1,...,5}
			\foreach \j in {1,...,5}
				\path[-] (I-\i) edge (N1-\j);

          \foreach \name / \y in {1}
			\path[yshift=(-3)*0.5cm]
				node[neuron, pin={[pin edge={->, dashed, black, very thick}]right:$\mathrm{Prob}(w_N=1)$}] (N2-\name) at (6, -1.6*\y cm) {
					$\mathsf{n}_\y^{(2)}$
         		};

		\foreach \i in {1,...,5}
			\foreach \j in {1}
				\path (N1-\i) edge (N2-\j);

		\node at (1.5+7.2, -1cm) {$\mathrm{ReLu}$};
		\draw (1.5+6,-1cm) -- (1.5+6.5,-1cm);

		\node at (1.5+6.9, -1.7cm) {$\sigma$};
		\draw[dashed,->] (1.5+6,-1.7cm) -- (1.5+6.5,-1.7cm);
	\end{tikzpicture}
	\caption{Neural network architecture used for predicting the Fricke sign from feature vectors of the form $(a_2,a_3,\dots,a_{p_d},R)$
    where $p_i$ denotes $i^{\mathrm{th}}$ prime number, and $R$ denotes the spectral parameter.
     Here $\mathsf{n}_{i}^{(j)}$ denotes $i^{\mathrm{th}}$ node in $j^{\mathrm{th}}$ layer.
    The spectral parameter $R$ is switched off for generating one of the visualizations in Figure~\ref{fig:acc_vs_params_NN}.}
	\label{fig:nnArchitecture}
\end{figure}
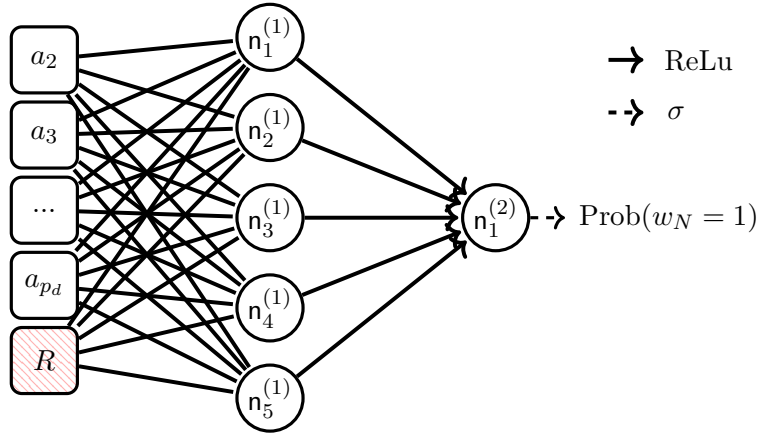

The neural network architecture is shown in Figure~\ref{fig:nnArchitecture}.
We train the networks over $4\times 10^{4}$ iterations, using Adam optimization with learning rate $10^{-3}$.
This achieves around $95\%$ accuracy for even forms and $94\%$ accuracy for odd forms.
For the loss function, we use binary cross-entropy, with the classes corresponding to two different possibilities for the Fricke sign $w_N$.

\begin{figure}[htb]
    \centering
    \includegraphics[width=0.8\textwidth]{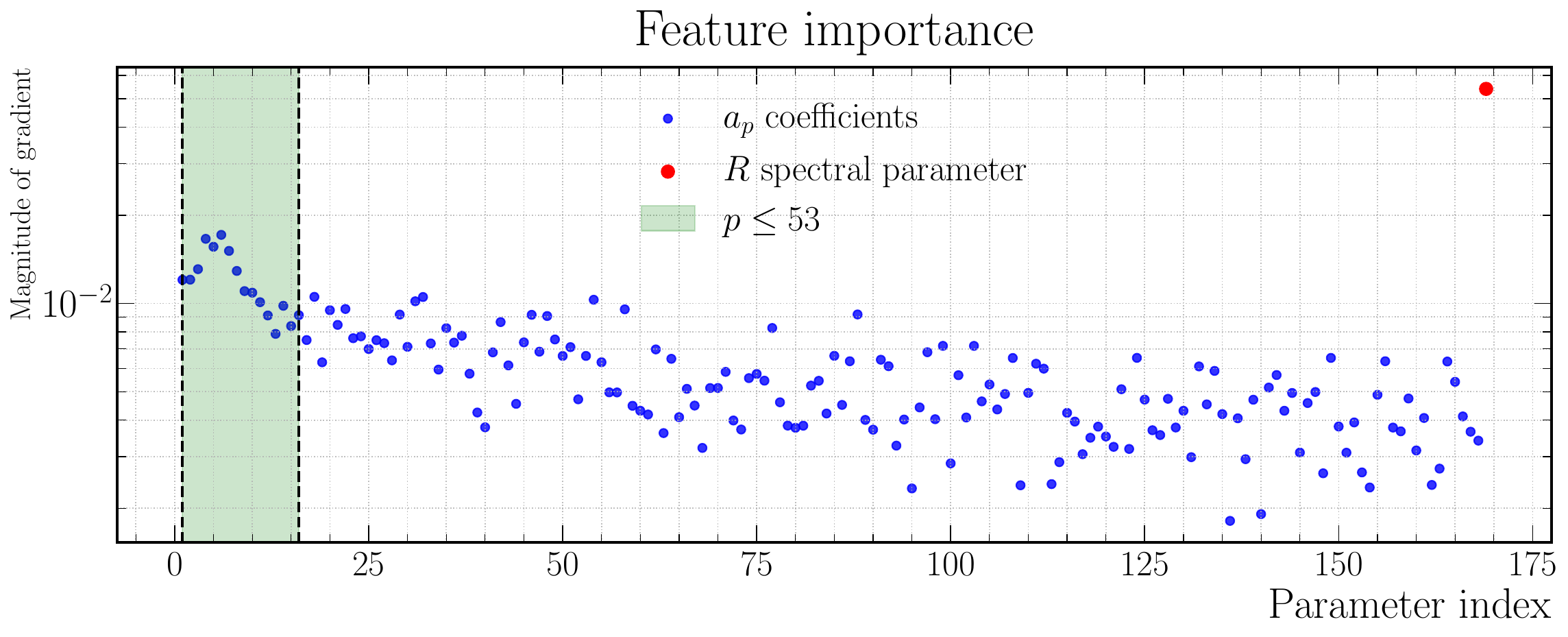}
    \caption{Saliency map of a neural network trained for predicting Fricke sign.}
    \label{fig:saliency_NN}
\end{figure}

\begin{figure}
	\centering
	\includegraphics[width=0.8\textwidth]{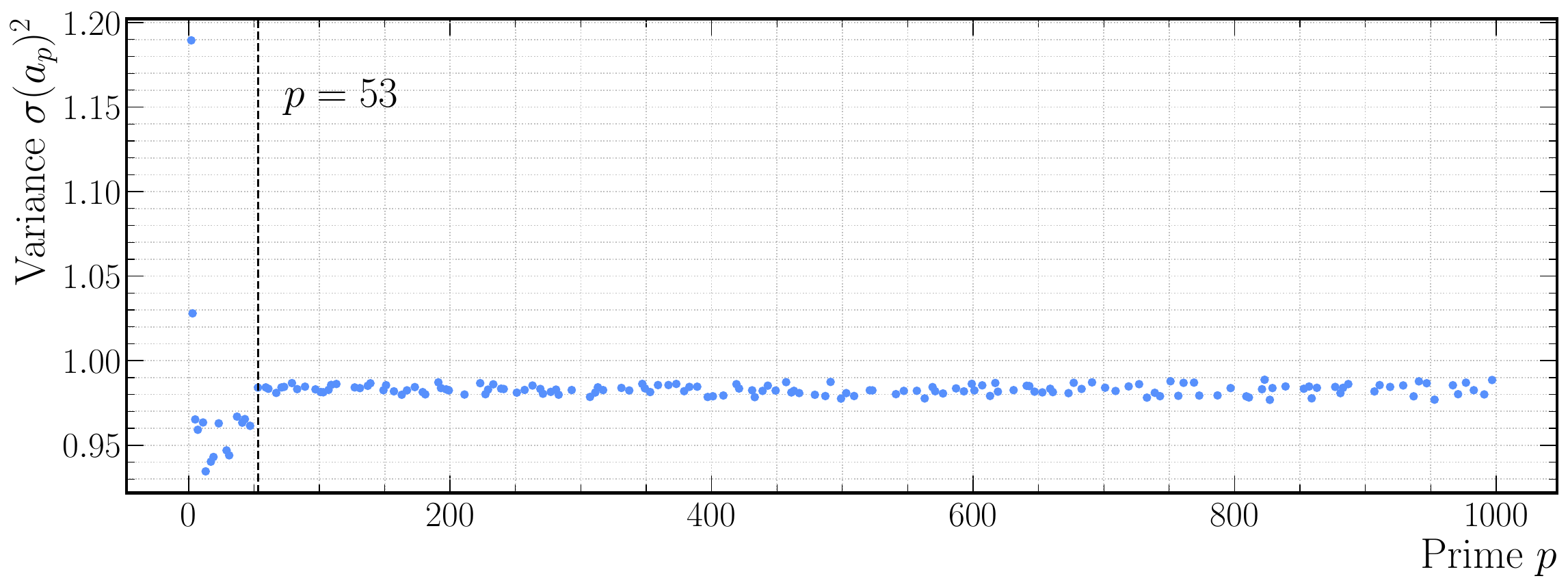}
	\caption{Variances of the coefficients $a_p$ for Maass forms from \cite{lmfdb}.
    }
	\label{fig:variances}
\end{figure}

The saliency analysis of the trained neural network (shown in Figure~\ref{fig:saliency_NN}) indicates the coefficients $\{a_p~\vert~ p \leq 53\}$ and the spectral parameter $R$ contribute the most to the decision of the neural network.
While the input parameters, including the coefficients $(a_{p_i})_{i=1}^d$ and the spectral parameter $R$, are normalized to have unit variance over the dataset, the standard deviations of the distributions of the coefficients indicate a sharp jump around $p=53$, as shown in Figure~\ref{fig:variances}.
This variance is caused by the
$a_p$ having only three distinct values in our dataset, $\{0, \pm 1/\sqrt{p}\}$,  when $p$ divides the level. This indicates that it is possible that the neural network architecture is making more use of the way the Fricke sign is embedded in the the Fourier coefficients when $p|N$.
The Maass forms in the LMFDB have analytic conductor up to about $600$.
Motivated by this observation, we measure the influence that the number of the coefficients has on the performance of the network.
In particular, we train each neural network on a subset of coefficients
of varying size, and calculate the resulting accuracy under the same choice of hyperparameters.
The dependence of the accuracy on the number of coefficients is shown on Figure~\ref{fig:acc_vs_params_NN}.

\subsection{Comparison of Hejhal's algorithm, LDA, and neural networks}\label{s:comparisonLMFDB}

The rigorous computations of Maass forms in the LMFDB were computationally expensive.
It is easier to make \emph{heuristic} guesses at the Fricke signs.
One approach is to modify (the original, heuristic, version of) Hejhal's algorithm~\cite{hejhal} to apply to general level.
This works by first guessing an eigenvalue and a set of Atkin--Lehner involution signs.
Then, one evaluates a Maass form $f$ with that eigenvalue and Atkin--Lehner
behavior at several points, and also evaluates $f$ at images of these points
under $f(z) \mapsto f(\gamma z)$ for $\gamma \in \Gamma_0(N)$.
If the guessed eigenvalue is close enough to be correct, and if the Atkin--Lehner signs are correct, then these expansions must approximately agree by the automorphy of $f$.
By truncating the Fourier expansions, this becomes an over-determined, approximate, homogenous linear system in the coefficients.
The non-linearity of the group action non-trivially mixes the coefficients to make the system of equations solvable in practice.

We ran this heuristic Hejhal's algorithm as given in~\cite{dldgithub} for each of the Maass forms in the LMFDB with unknown Fricke sign.
As the Fricke signs were unknown, it was necessary to try every possible combination of Atkin--Lehner signs, which makes the computation much more expensive.
For each eigenvalue and set of signs, Hejhal's algorithm produces a candidate list of coefficients.
We can heuristically guarantee that these coefficients are close to correct if they satisfy Hecke relations and are close to known coefficients.

This heuristic algorithm yields expansions and Fricke signs that are heuristically extremely close to correct for 4,595 of the 15,423 Maass forms in the LMFDB with unknown Fricke signs.
We consider these as ``probably correct'' and compare the models here against this data.

\begin{remark}
    Hejhal's intention was to iterate the algorithm to find spectral eigenvalues and their Maass forms.
    For each choice of Atkin--Lehner signs, we can regard Hejhal's algorithm as a demanding root-finding algorithm in the unknown spectral eigenvalue.
    For those forms with unknown Fricke sign, it would be necessary to
    repeatedly run this algorithm with every combination of Atkin--Lehner signs
    and with several eigenvalue candidates.
    It was surprising to the authors that the current precision in the LMFDB isn't sufficient for even heuristic versions of Hejhal's algorithm to converge in practice.
\end{remark}

\begin{table}[htb!]
\begin{tabular}{|c|c|}
\hline
LDA features & $\%$ agreement with Hejhal heuristic  \\
\hline
$a_n$ &  95.45\\
\hline
$a'_n$  & 94.84\\
\hline
$a_{p_i}$ & 82.93\\
\hline
$a'_{p_i}$ & 82.87 \\
\hline
\end{tabular}
\caption{Comparing the LDA predictions with the heuristic approach based on Hejhal's algorithm. }
\label{tab:LDA-heuristic}
\end{table}

\begin{table}[ht!]
\begin{tabular}{|c|c|}
\hline
Parity & $\%$ agreement with Hejhal heuristic  \\
\hline
odd  &  96.09  (1877 out of 1984)\\
\hline
even &  94.61 (2509 out of 2611)\\
\hline
\end{tabular}
\caption{Comparison of LDA with the heuristic approach based on Hejhal's algorithm for Maaass forms with fixed parity. For LDA, we used feature vectors of the form $(a_n)_{n=1}^{1000}$. }
\label{tab:LDA-heuristic-fixed-par}
\end{table}

\begin{table}[htb!]
\begin{tabular}{|c|c|}
\hline
Parity & $\%$ agreement with Hejhal heuristic  \\
\hline
both &  $87.37$\\
\hline
odd &  $85.79$\\
\hline
even & $89.46$\\
\hline
\end{tabular}
\caption{Comparing the neural network predictions with the heuristic approach based on Hejhal's algorithm.}
\label{tab:NN-heuristic}
\end{table}

In Table~\ref{tab:LDA-heuristic}, we compare the predictions using the 4 different LDA methods from Table~\ref{tab:no_cheat} with those from the heuristic Hejhal algorithm.
The features that achieve the greatest agreement are the standard Fourier coefficients $a_n$.
The performance on odd Maass forms is slightly better than on even ones, similar to results on the validation set.
The high accuracy is sustained when we restrict to forms of fixed parity.
For example, using $a_n$ as features, we achieve $95.45\%$ accuracy on all forms, and, if one only looks at odd (resp. even) forms, then the accuracy is $96.09\%$ (resp. $94.61\%$).
A comparison between LDA and Hejhal's heuristic for Maass forms with fixed parity is documented in Table~\ref{tab:LDA-heuristic-fixed-par}.
In Table~\ref{tab:NN-heuristic}, we compare the predictions using neural networks with those from the heuristic Hejhal's algorithm.

\bibliographystyle{alpha}
\bibliography{bibfile}

\end{document}